\documentclass[conference]{IEEEtran}
\IEEEoverridecommandlockouts

\usepackage[T1]{fontenc}
\usepackage[latin9]{inputenc}
\usepackage[numbers, compress]{natbib}
\usepackage{amsthm}
\usepackage{amsmath}
\usepackage{algcompatible}
\usepackage{algorithm}
\usepackage{booktabs}
\usepackage{enumerate}
\usepackage{graphicx}
\usepackage{amssymb}
\usepackage{latexsym}
\usepackage{epstopdf}
\usepackage{color}
\usepackage{xcolor}
\usepackage{bbm}
\usepackage{needspace}
\usepackage{color, colortbl}
\usepackage[english]{babel}
\usepackage{tikz}
\usepackage{caption}
\usepackage{subcaption}
\usetikzlibrary{arrows,automata}
\usetikzlibrary{positioning}
\usepackage{filecontents}
\usepackage{mathtools}
\DeclarePairedDelimiter{\ceil}{\lceil}{\rceil}

\DeclareCaptionFont{mysize}{\fontsize{8}{9.6}\selectfont}
\ifCLASSINFOpdf
\else
\fi

\theoremstyle{plain}
\newtheorem{thm}{\protect\theoremname}
  \theoremstyle{plain}

\newtheorem{define}{Definition}
\newtheorem{problem}{Problem}

\makeatother

\usepackage{babel}
\providecommand{\lemmaname}{Lemma}
\providecommand{\theoremname}{Theorem}

\newcommand{\bdA}{\ensuremath{\boldsymbol{\alpha}}}

\begin{document}
\title{\vspace{0.1in} Safety Verification of Nonlinear Autonomous System via Occupation Measures}
\author{Ximing Chen, Shaoru Chen, Victor M. Preciado%
\thanks{The authors are with the Department of Electrical and Systems Engineering at the University of Pennsylvania, Philadelphia PA 19104. e-mail:\{ximingch, srchen, preciado\}@seas.upenn.edu. %
}}

\maketitle

\begin{abstract}
In this paper, we introduce a flexible notion of safety verification for nonlinear autonomous systems by measuring how much time the system spends in given unsafe regions. We consider this problem in the particular case of nonlinear systems with a polynomial dynamics and unsafe regions described by a collection of polynomial inequalities. In this context, we can quantify the amount of time spent in the unsafe regions as the solution to an infinite-dimensional linear program (LP). This LP measures the volume of the unsafe region with respect to the occupation measure of the system trajectories. Using Lasserre hierarchy, we approximate the solution to the infinite-dimensional LP using a sequence of finite-dimensional semidefinite programs (SDPs). The solutions to the SDPs in this hierarchy provide monotonically converging upper bounds on the optimal solution to the infinite-dimensional LP. Finally, we validate the performance of our framework using numerical simulations.
\end{abstract}

\section{Introduction}\label{sec:introduction}
Our ability to provide safety certificates about the behavior of complex systems is critical in many engineering applications, such as air traffic control~\cite{hu2003probabilistic}, life support devices~\cite{glavaski2005safety}, motion planning in robotics manipulations~\cite{ziegler2010fast}, and connected autonomous vehicles~\cite{althoff2007safety,althoff2009model}. Although safety verification is a mature area with many success stories~\cite{bemporad2000optimization,chutinan2003computational}, the verification of nonlinear dynamical systems over nonconvex unsafe regions remains a challenging problem~\cite{bollobas1997volume,dyer1988complexity}. 

In the past decades, various solutions have been proposed to verify the safety of dynamical systems. The solution approaches often fall into the following two categories: (i) reachable set methods~\cite{anai2001reach, tomlin2003computational, asarin2003reachability}, and (ii) Lyapunov function methods~\cite{prajna2004safety,prajna2006barrier,prajna2007framework,sloth2012compositional}. Essentially, reachable set methods aim to find a set containing all possible states at a given time, for a given set of initial conditions. Subsequently, if the reachable set does not intersect with the pre-specified unsafe regions, the system is considered to be safe. For example, in~\cite{anai2001reach} the reachable set is found for continuous-time linear systems, whereas in~\cite{tomlin2003computational} and \cite{asarin2003reachability} the reachable sets are computed via approximations for nonlinear dynamical systems. In~\cite{kousik2017safe}, the authors applied a reachable set method to plan safe trajectories for autonomous vehicles.

While reachable set methods can be used to obtain quantitative guarantees for safety, the reliability of the result largely depends on the assumptions made about the system, as well as the form of the unsafe regions. For instance, calculating the volume of the intersection of two sets, such as the reachable set and the unsafe regions, can become computationally challenging~\cite{dyer1988complexity}, jeopardizing the practical application of reachable set methods. An alternative approach to safety verification is based on using Lyapunov-like functions. In~\cite{prajna2006barrier}, the authors proposed the use of barrier certificates for safety verification of nonlinear systems. In contrast with the reachable set method, this line of work does not require to solve differential equations and is computationally more tractable. Furthermore, it also allows to provide safety certificates for various types of hybrid~\cite{prajna2004safety} and stochastic systems~\cite{prajna2007framework}. 

Despite a tremendous amount of solutions proposed to solve the safety verification problem, the majority of existing methods only provide binary safety certificates. More specifically, these certificates concern only \emph{whether the system is safe} rather than \emph{how safe the system is}. Lacking a detailed analysis of how unsafe a system is may result in a restricted and conservative design space. To illustrate this point, let us consider the operation of a solar-powered autonomous vehicle. Naturally, regions without solar exposure are considered to be unsafe, since the battery of the vehicle could be drained after a period of time. However, it would be inefficient to plan a path for the vehicle completely avoiding all these shaded regions. Instead, a more suitable requirement would be that the amount of time the vehicle spends in the shaded regions is bounded. More generally, this framework can be useful in those situations where the system is able to tolerate the exposure to a deteriorating agent, such as excessive heat or radiation, for a limited amount of time.

In this paper, we consider this alternative, more flexible notion of safety. More precisely, we aim to compute the time that a (nonlinear) system spends in the unsafe regions. In particular, we focus our analysis on the case of systems described by a polynomial dynamics and unsafe regions described by a collection of polynomial inequalities. To calculate the amount of time spent in the unsafe regions, we use \emph{occupation measures} to quantify how much time the system trajectory spends in a particular set~\cite{vinter1993convex}. Using this alternative viewpoint of the system dynamics, the safety quantity of interest can be calculated by finding the volume of the unsafe region with respect to the occupation measure~\cite{henrion2009approximate}. The usage of occupation measures allows us to leverage powerful numerical procedures developed in the context of control of polynomial systems~\cite{lasserre2008nonlinear,henrion2014convex,majumdar2014convex}.

The contribution of this paper is threefold. First, we formulate a flexible notion of safety allowing a trade-off between safety and performance. Second, we provide an \emph{exact} formulation of the problem under consideration in terms of an infinite-dimensional LP. Furthermore, we provide a hierarchy of relaxations that can be efficiently solved using semidefinite programming. Finally, we provide numerical examples to demonstrate the applicability of our method. 

The rest of the paper is structured as follows. The safety verification problem formulation is stated in Section II. In Section III, we introduce concepts from measure theory that are necessary for developing our framework. Based on those notions, we show that the problem under consideration can be stated as an infinite-dimensional linear program, and in Section IV, we provide approximate solutions to this LP using a sequence of semidefinite programs. The performance of our framework is illustrated through numerical experiments in Section V and we conclude the paper in Section VI.

\textbf{Notations:} We use bold symbols to represent real-valued vectors. Given $n\in \mathbb{N},$ we use the shorthand notation $[n]$ to denote the set of integers $\{1,\ldots, n\}.$ The indicator function of a given set $\mathcal{S}$ is defined by $\mathbf{1}_\mathcal{S}(\cdot).$ We use $\delta_\mathbf{x}$ to denote the Dirac measure centered on a fixed point $\mathbf{x} \in \mathbb{R}^n$ and we use $\otimes$ to denote the product between two measures. The ring of polynomials in $\mathbf{x}$ with real coefficients is denoted by $\mathbb{R}[\mathbf{x}]$, and $\mathbb{R}[\mathbf{x}]_r$ denotes the subset of polynomials of degree $\leq r$. Given $\mathbf{x}\in\mathbb{R}^n$ and $\bdA \in \mathbb{N}^n,$ we let $\mathbf{x}^{\boldsymbol{\alpha}}$ denote the quantity $\mathbf{x}^{\bdA} = \prod_{i = 1}^n x_i^{\alpha_i}$. Let $|\boldsymbol{\alpha}| = \sum_{i = 1}^n \alpha_i$ and $\mathbb{N}_r^n = \lbrace \bdA \in \mathbb{N}^n \mid |\bdA| \leq r \rbrace$.

\section{Problem Statement}\label{sec:statement}
In this paper, we consider a continuous-time autonomous dynamical system whose dynamics is captured by the following equation:
\begin{equation}\label{eq:dynamics}
\begin{aligned}
& \dot{\mathbf{x}}(t) = f(t, \mathbf{x}), \quad t\in [0, T]\\
& \mathbf{x}(0) = \mathbf{x}_0 \\
\end{aligned}
\end{equation}
where $\mathbf{x}(t) \in \mathbb{R}^n$ is the state vector, $\mathbf{x}_0$ is the initial condition, and $T > 0$ is the terminal time. We consider that the states of \eqref{eq:dynamics} are constrained to live within the set $\mathcal{X} \subseteq \mathbb{R}^n$ for all  $t \in [0,T].$ Furthermore, we consider that the system evolves from an initial condition $\mathbf{x}_0$, with $\mathbf{x}_0 \in \mathcal{X}_0 \subseteq \mathcal{X}.$ In this paper, we are interested in the case that the set $\mathcal{X}$ is semi-algebraic, as stated below. 

\begin{define}
A set $K\subseteq \mathbb{R}^n$ is said to be \emph{semi-algebraic} if there exist $m$ polynomials, $g_i: \mathbb{R}^n \rightarrow \mathbb{R}$, such that 
\begin{equation}\label{eq:SemiAlgebraicSet}
K = \{\mathbf{x}\in \mathbb{R}^n \mid g_i(\mathbf{x}) \ge 0, \forall i\in [m]\}.
\end{equation}
\end{define} 
As mentioned above, we assume that the set $\mathcal{X}$ can be defined using polynomials $g_i^{\mathcal{X}}(\mathbf{x}) \in \mathbb{R}[\mathbf{x}]$, as follows:
\begin{equation}\label{eq:setconstraint}
\begin{aligned}
&\mathbf{x}(t) \in \mathcal{X} = \{ \mathbf{x}\in \mathbb{R}^n \mid g_i^{\mathcal{X}}(\mathbf{x}) \ge 0, \forall i \in [n_{\mathcal{X}}]\} \\
\end{aligned}
\end{equation}
for all $t\in [0, T].$ 

In this paper, we consider the following problem:
\begin{problem}\label{Prob:Safety}	
Consider a \emph{compact} and \emph{semi-algebraic} set $\mathcal{X}$, defined by \eqref{eq:setconstraint}, and $\mathcal{X}_u \subseteq \mathcal{X}$, defined by:
\begin{equation}
\mathcal{X}_u = \{ \mathbf{x}\in \mathbb{R}^n \mid g_i^{\mathcal{X}_u}(\mathbf{x}) \ge 0, \forall i \in [n_{\mathcal{X}_u}]\}.
\end{equation}
Given the autonomous system described in~\eqref{eq:dynamics}, with $x_0 \sim \mu_0(\mathcal{X}_0)$, where $\mu_0$ is a probability distribution supported on $\mathcal{X}_0$, compute the expected amount of time that the system trajectory spends in the unsafe region $\mathcal{X}_u$. 
\end{problem}
Notice that this expected time can be computed as:
\begin{equation}\label{eq:probsafe0}
\mathbb{E}\left[\int_{0}^T \mathbf{1}_{\mathcal{X}_u}(\mathbf{x}(t)) dt\right],
\end{equation}
where the expectation in \eqref{eq:probsafe0} is taken with respect to the distribution of the initial condition $\mathbf{x}_0$. We remark that the above formulation is also capable of providing safety certificate for the system when the initial state is known exactly, i.e., $\mu_0 = \delta_{\mathbf{x}_0}.$

\section{Occupation Measure-based Reformulation}\label{sec:prelim}
In this section, we introduce a measure-theoretic approach to characterize the trajectories of the autonomous system described in~\eqref{eq:dynamics} presented in Subsection~\ref{subsec:OM}. Using this method, we show that the expectation in~\eqref{eq:probsafe0} can be computed via an infinite-dimensional linear program -- see Subsection~\ref{subsec:IDLP} and Subsection~\ref{subsec:DIDLP}. To explain our approach, we first introduce some notions of measure theory.

\subsection{Notations and preliminaries}\label{subsec:notations}
Given a topological space $\mathcal{S},$ we denote by $\mathcal{M(S)}$ the space of finite signed Borel measures on $\mathcal{S}$, and $\mathcal{M}_{+}(\mathcal{S})$ its positive cone. Let $\mathcal{C(S)}$ and $\mathcal{C}^1(\mathcal{S})$ be the space of continuous functions and continuously differentiable functions on $\mathcal{S}$, respectively. The topological dual of $\mathcal{M(S)}$ and $\mathcal{C(S)}$ are denoted by $\mathcal{M(S)}^*$ and $\mathcal{C(S)}^*$.

Given a function $h \in \mathcal{C}(\mathcal{S})$ and a measure $\mu \in \mathcal{M}(\mathcal{S}),$ we define the duality bracket between $h$ and $\mu$ by
\begin{equation}\label{eq:dualitybracket}
\langle h, \mu \rangle = \int_{\mathcal{S}} hd\mu.
\end{equation}
By Riesz-Markov-Kakutani representation theorem~\cite{kakutani1941concrete}, when $\mathcal{S}$ is locally compact Hausdorff, the dual space of $\mathcal{C}(\mathcal{S})$ is $\mathcal{M(S)}$, in which the norm of $\mathcal{C}(\mathcal{S})$ is the $\sup$-norm of functions and the norm of $\mathcal{M}(\mathcal{S})$ is the total variation norm of measures. In the rest of the paper, we consider compact topological spaces $\mathcal{S} \subseteq \mathbb{R}^n$. As a consequence, both local compactness and separability conditions required to form the duality between $\mathcal{M(S)}$ and $\mathcal{C(S)}$ are satisfied. Given a measure $\mu \in \mathcal{M}(\mathcal{S}),$ the support of $\mu,$ denoted by $\mathtt{supp}(\mu)$, is the smallest closed set $C \subseteq \mathcal{S}$ such that $\mu(\mathcal{S} \setminus C) = 0$ where smallest is understood in the set-inclusion sense.

\subsection{Occupation measures and Liouville equation}\label{subsec:OM}
Given an initial condition $\mathbf{x}_0,$ let $\mathbf{x}(t \mid \mathbf{x}_0)$ be the solution to~\eqref{eq:dynamics}. Given a trajectory $\mathbf{x}(t \mid \mathbf{x}_0),$ we define the \emph{occupation measure} $\mu(\cdot \mid \mathbf{x}_0)$ of $\mathbf{x}(t \mid \mathbf{x}_0)$ as 
\begin{equation}\label{eq:OccupationMeasure}
\mu(A \times B \mid \mathbf{x}_0) = \int_{[0, T]\cap A} \mathbf{1}_{B} (\mathbf{x}(t \mid \mathbf{x}_0)) dt
\end{equation}
for all $A \times B \subseteq [0, T]\times \mathcal{X}.$ Therefore, given sets $A$ and $B,$ the value $\mu(A \times B)$ equals the total amount of time out of $A$ that the state trajectory $\mathbf{x}(t \mid \mathbf{x}_0)$ spends in the set $B$. 
Similarly, we define the \emph{final measure} $\mu_T(\cdot \mid\mathbf{x}_0)$ as
\begin{equation}\label{eq:FinalMeasure}
\mu_T(B \mid \mathbf{x}_0) = \mathbf{1}_B(\mathbf{x}(T \mid \mathbf{x}_0))
\end{equation}
for $B \subseteq \mathcal{X}.$ Notice that the occupation measure $\mu(\cdot \mid \mathbf{x}_0)$ is supported on $[0,T] \times \mathcal{X}$ whereas the final measure $\mu_T(\cdot \mid \mathbf{x}_0)$ is supported on $\mathcal{X}.$ 

Given a test function $v \in \mathcal{C}^1([0, T] \times \mathcal{X}),$ we define the operator $\mathcal{L}$ as:
\begin{equation}\label{eq:Operator}
v \mapsto \mathcal{L}v = \dfrac{\partial v}{\partial t} +  \nabla v \cdot f(t, \mathbf{x}).
\end{equation}
The \emph{adjoint operator} $\mathcal{L}^*: \mathcal{M}([0, T] \times \mathcal{X}) \rightarrow \mathcal{C}^1([0, T] \times \mathcal{X})^*$ is given by
\begin{equation}\label{eq:AdjointOperator}
\langle v, \mathcal{L}^* \nu \rangle = \langle \mathcal{L}v, \nu \rangle.
\end{equation}
From~\eqref{eq:Operator}, we have that
\begin{equation}\label{eq:derivation}
\begin{aligned}
v(T, \mathbf{x}(T \mid \mathbf{x}_0)) & = v(0, \mathbf{x}_0) + \int_{0}^T \dfrac{d}{dt}v(t, \mathbf{x}(t \mid \mathbf{x}_0)) dt \\
&= v(0, \mathbf{x}_0) + \int_{[0, T]\times \mathcal{X}} \mathcal{L}v(t, \mathbf{x})d\mu(t, \mathbf{x} \mid \mathbf{x}_0) \\
&= v(0, \mathbf{x}_0) + \langle \mathcal{L} v, \mu(\cdot\mid \mathbf{x}_0)\rangle.
\end{aligned}
\end{equation}
Hence, we can further rewrite~\eqref{eq:derivation} as
\begin{equation}
\langle v, \delta_{T} \otimes \mu_T(\cdot \mid \mathbf{x}_0)  \rangle = \langle v, \delta_{0} \otimes \delta_{\mathbf{x}_0} \rangle + \langle \mathcal{L}v, \mu(\cdot\mid \mathbf{x}_0)\rangle.
\end{equation}
In the view of~\eqref{eq:AdjointOperator}, since the above equation holds for all $v \in \mathcal{C}^1([0, T]\times \mathcal{X}),$ we obtain the following equality:
\begin{equation}\label{eq:Liouville}
\delta_{T} \otimes \mu_T(\cdot \mid \mathbf{x}_0) = \delta_{0} \otimes \delta_{\mathbf{x}_0} + \mathcal{L}^* \mu(\cdot\mid \mathbf{x}_0).
\end{equation}
Essentially,~\eqref{eq:Liouville} describes the evolution of the distribution of states, given an initial distribution, under the flow of the dynamics~\eqref{eq:dynamics} -- see~\cite{arnol2013mathematical} for a more detailed discussions.

The measures defined in~\eqref{eq:OccupationMeasure} and~\eqref{eq:FinalMeasure} depend on a given initial condition $\mathbf{x}_0.$ In what follows, we extend these definitions to handle the case when the system is evolving from a set of possible initial conditions. Given an initial distribution $\mu_0$ with $\mathtt{supp}(\mu_0) \subseteq \mathcal{X}_0$, we define the \emph{average occupation measure} $\mu \in \mathcal{M}([0, T] \times \mathcal{X})$ as
\begin{equation}\label{eq:AverageOccupationMeasure}
\mu(A \times B) = \int_{\mathcal{X}_0} \mu(A \times B \mid \mathbf{x}_0) d\mu_0
\end{equation}
and the \emph{average final measure} $\mu_T \in \mathcal{M}(\mathcal{X})$ as
\begin{equation}\label{eq:AverageFinalMeasure}
\mu_T(B) = \int_{\mathcal{X}_0} \mu_T(B\mid \mathbf{x}_0)d\mu_0.
\end{equation}
By integrating the left- and right-hand side of~\eqref{eq:derivation} with respect to $\mu_0$, we have that
\begin{equation}\label{eq:AverageLiouvilleLinear}
\delta_T \otimes \mu_T = \delta_0 \otimes \mu_0 + \mathcal{L}^*\mu.
\end{equation}

Note that any family of solutions $\mathbf{x}(t)$ of \eqref{eq:dynamics} with an initial distribution $\mu_0$ induces an occupation measure \eqref{eq:AverageOccupationMeasure} and a final measure \eqref{eq:AverageFinalMeasure} satisfying \eqref{eq:AverageLiouvilleLinear}. Conversely, for any tuple of measures $(\mu_0, \mu, \mu_T)$ satisfying~\eqref{eq:AverageLiouvilleLinear}, one can identify a distribution on the admissible trajectories starting from $\mu_0$ whose average occupation measure and average final measure coincide with $\mu$ and $\mu_T$, respectively (see Lemma~$3$ in \cite{henrion2014convex} and Lemma~$6$ in \cite{zhao2017control} for more details). 

\subsection{Infinite-dimensional linear program reformulation}
\label{subsec:IDLP}
Hereafter, we will show that the value in~\eqref{eq:probsafe0} can be obtained by solving a linear program on the occupation measure and the final measure, defined in~\eqref{eq:AverageOccupationMeasure} and~\eqref{eq:AverageFinalMeasure}. According to the definition of average occupation measure, we have that
\begin{equation}
\begin{aligned}
\mathbb{E}\left[\int_{0}^T \mathbf{1}_{\mathcal{X}_u}(\mathbf{x}(t)) dt\right] &= \int_{\mathcal{X}_0}\int_{0}^T \mathbf{1}_{\mathcal{X}_u}(\mathbf{x}(t)) dtd\mu_0 \\
&= \int_{\mathcal{X}_0}\mu([0, T] \times \mathcal{X}_u \mid \mathbf{x}_0) d\mu_0 \\
&= \mu([0, T]\times \mathcal{X}_u).
\end{aligned}
\end{equation}
Leveraging the above measure-theoretical formulation, the value in~\eqref{eq:probsafe0} is equal to
\begin{equation}\label{eq:probsafe}
\mu([0,T]\times\mathcal{X}_u).
\end{equation}
Subsequently, finding the solution to Problem~\ref{Prob:Safety} is equivalent to finding the \emph{volume} of the set $[0, T]\times \mathcal{X}_u$, where this volume is measured using the average occupation measure, instead of the Lebesgue measure. Next, we show that the value of~\eqref{eq:probsafe} can be obtained by solving the following optimization problem: Given a polynomial $g: [0,T]\times \mathcal{X} \rightarrow \mathbb{R},$ such that $g(t, \mathbf{x}) > 0, \forall (t, \mathbf{x}) \in [0,T] \times \mathcal{X}_u,$ consider the following optimization problem
\begin{equation}\label{eq:MainProgram}
\begin{aligned}
\mathtt{P}: \indent \indent \indent & \sup \int g d\widetilde{\mu} \\
& \text{subject to } \widetilde{\mu} + \widehat{\mu} = \mu\\
&\indent \indent \indent \indent \delta_T \otimes \mu_T = \delta_0 \otimes \mu_0 + \mathcal{L}^*\mu\\
&\indent \indent \indent \indent \mu, \widehat{\mu} \in \mathcal{M}_{+}([0,T]\times \mathcal{X}) \\
& \indent \indent \indent \indent \widetilde{\mu} \in \mathcal{M}_{+}([0,T]\times \mathcal{X}_u)\\
& \indent \indent \indent \indent \mu_T \in  \mathcal{M}_{+}(\mathcal{X})\\
\end{aligned}
\end{equation}
where the supremum is taken over a tuple of measures $(\widetilde{\mu}, \widehat{\mu}, \mu, \mu_T)  \in \mathcal{M}_{+}([0,T] \times \mathcal{X}_u) \times \mathcal{M}_{+}([0,T] \times \mathcal{X}) \times \mathcal{M}_{+}([0,T] \times \mathcal{X}) \times \mathcal{M}_{+}(\mathcal{X}).$ The constraint $\widetilde{\mu} + \widehat{\mu} = \mu$ is equivalent to $\widetilde{\mu} \le \mu$, i.e., the measure $\widetilde{\mu}$ is dominated by $\mu$. Using duality brackets, we can write the objective in~\eqref{eq:MainProgram} as $\langle g, \widetilde{\mu} \rangle.$ It follows that~\eqref{eq:MainProgram} is a linear program in the decision variable $(\widetilde{\mu}, \widehat{\mu}, \mu, \mu_T)$. Denote by $\sup \mathtt{P}$ the optimal value of $\mathtt{P}$ and by $\max \mathtt{P}$ the supremum attained. When $g \equiv 1$, we show below that the optimal value to the above program, if it exists, is equal to~\eqref{eq:probsafe0}.

\begin{thm}\label{thm:LPOptimal}
Let $\mathcal{X}_u$ be a compact and semi-algebraic subset of $\mathcal{X}$ and $\mathcal{B}$ be the Borel $\sigma$-algebra of Borel subsets of $[0,T]\times\mathcal{X}.$ Let $\widetilde{\mu}^* \in \mathcal{M}([0, T] \times \mathcal{X}_u)$ be defined by
\begin{equation} \label{eq:VolumeDefinition}
\widetilde{\mu}^*(S) = \mu(S \cap [0, T] \times \mathcal{X}_u), \forall S \in \mathcal{B}.
\end{equation}
Given a polynomial $g: [0,T]\times \mathcal{X} \rightarrow \mathbb{R},$ if $g(t, \mathbf{x}) > 0, \forall (t, \mathbf{x}) \in [0, T] \times \mathcal{X}_u,$ then $\widetilde{\mu}^*$ is the $\widetilde{\mu}$-component of an optimal solution to $\mathtt{P}$. Furthermore, $\sup \mathtt{P} = \max \mathtt{P} = \int gd \widetilde{\mu}^*.$ In particular, if $g \equiv 1,$ then $\max \mathtt{P} = \mu([0, T]\times\mathcal{X}_u).$
\end{thm}
\begin{proof}
See appendix~\ref{appen:proof}.
\end{proof}

As a result of Theorem~\ref{thm:LPOptimal}, the solution of $\mathtt{P}$ is equal to the expected time in \eqref{eq:probsafe0}. In the next subsection, we consider the Lagrangian dual of $\mathtt{P}$.

\subsection{Dual infinite-dimensional program}\label{subsec:DIDLP}
As mentioned in Section~\ref{subsec:notations}, the dual space of $\mathcal{M}(\mathcal{S})$ is the Banach space of continuous functions on $\mathcal{S}$ with the sup-norm. Let $\mathcal{C}_+(\mathcal{S}) \subseteq \mathcal{C}(\mathcal{S})$ be the set of continuous functions that are nonnegative on $\mathcal{S}.$ Using duality theory, the dual program of~\eqref{eq:MainProgram} is equal to
\begin{equation}\label{eq:DualProgram}
\begin{aligned}
\mathtt{D}: \ &\underset{v,w}{\inf} \int v(0, \mathbf{x}) d\mu_0 \\
& \text{s.t.} \ w(t, \mathbf{x}) - g(t, \mathbf{x}) \geq 0, \forall (t, \mathbf{x}) \in [0, T]\times \mathcal{X}_u\\
&\indent \indent -\mathcal{L}v(t, \mathbf{x}) - w(t, \mathbf{x}) \ge 0, \forall (t, \mathbf{x}) \in [0, T]\times \mathcal{X}\\
&\indent \indent  v(T, \mathbf{x}) \ge 0, \forall \mathbf{x} \in \mathcal{X} \\
&\indent \indent w(t,\mathbf{x}) \geq 0, \forall (t, \mathbf{x}) \in [0,T] \times \mathcal{X}
\end{aligned}
\end{equation}
where the decision variables in the above program are the continuously differentiable function $v(t,\mathbf{x}) \in \mathcal{C}^{1}([0, T]\times \mathcal{X})$ and the continuous function $w(t, \mathbf{x}) \in \mathcal{C}([0, T]\times \mathcal{X})$. The dual problem $\mathtt{D}$ always provides an upper bound on the optimal value of the primal $\mathtt{P}$. In the sequel, we show that the optimal values of~\eqref{eq:MainProgram} and~\eqref{eq:DualProgram} are actually equal. Thus, \emph{strong duality} holds in this infinite-dimensional linear program.


\begin{thm}\label{thm:StrongDuality}
Let $p^\star$ and $d^\star$ be the optimal values of $\mathtt{P}$ and $\mathtt{D},$ respectively. Then, $p^\star = d^\star$, i.e., there is no duality gap between $\mathtt{P}$ and $\mathtt{D}$.
\end{thm}
\begin{proof}
See appendix~\ref{appen:proof}.
\end{proof}

Consequently, the value of~\eqref{eq:probsafe0} can be obtained by solving~\eqref{eq:MainProgram} or~\eqref{eq:DualProgram}. However, these two optimization problems are taking arguments from a tuple of measures or a tuple of continuous functions; hence both programs are hard infinite-dimensional optimization problems. In the next section, we leverage recent results from the multi-dimensional moment problem~\cite{bernard2009moments} to approximate the solution to~\eqref{eq:MainProgram}. Furthermore, we show that it is possible to obtain increasingly tighter bounds on~\eqref{eq:probsafe} by solving a sequence of semidefinite programs.

\section{Semidefinite and Sum-of-Squares Relaxation}\label{sec:relaxation}
In the previous section, we have shown that~\eqref{eq:probsafe0} can be computed by solving an infinite-dimensional linear program. Although the optimal solutions to $\mathtt{P}$ or $\mathtt{D}$ provide exact solutions to Problem~1, it is computationally intractable to solve them. To address this issue, in Subsection~\ref{subsec:approximate}, we will provide a method to approximate the optimal solutions to $\mathtt{P}$ and $\mathtt{D}$ using sequences of semidefinite programs (SDPs) and sum-of-squares (SOS) programs, respectively. We utilize tools developed in the context of the multi-dimensional moment problem allowing us to replace the tuple of measures in $\mathtt{P}$ by sequences of moments. 

The following observation plays a key role in our approximation scheme. 
Notice that the equality constraint in~\eqref{eq:MainProgram} is equivalent to
\begin{equation}\label{eq:observation}
\langle v, \delta_T \otimes \mu_T\rangle  = \langle v, \delta_0 \otimes \mu_0\rangle + \langle \mathcal{L}v, \mu \rangle
\end{equation}
for all $v \in \mathcal{C}([0, T]\times \mathcal{X}).$ Since the set of polynomials are dense in $\mathcal{C}([0, T] \times \mathcal{X})$ and the ring $\mathbb{R}[t, \mathbf{x}]$ is closed under addition and multiplication,~\eqref{eq:observation} is equivalent to
\begin{equation}\label{eq:LiouvillePolynomial}
\begin{aligned}
& \int_{\mathcal{X}} v(T, \mathbf{x})d\mu_T = \int_{\mathcal{X}} v(0, \mathbf{x})d\mu_0 + \int_{[0, T]\times \mathcal{X}} \mathcal{L}v d\mu\\
& \text{ for all } v(t, \mathbf{x}) = t^a\mathbf{x}^{\bdA}, (a, \bdA) \in \mathbb{N}\times\mathbb{N}^n,
\end{aligned}
\end{equation}
where $a \in \mathbb{N}$, $\bdA = (\alpha_1, \cdots, \alpha_n)\in \mathbb{N}^n$ and $\mathbf{x}^{\bdA} = x_1^{\alpha_1}x_2^{\alpha_2}\cdots x_n^{\alpha_n}$. Using the above procedure, the linear constraints in $\mathtt{P}$ hold provided that~\eqref{eq:observation} holds for all monomial functions $v(t,\mathbf{x}).$ A standard relaxation is then to require that~\eqref{eq:observation} holds for all monomials up to a given fixed degree $r,$ i.e., $ a + |\bdA| = a + \sum_{i = 1}^n \alpha_i \le r.$ 

Since $v(t,\mathbf{x})$ is a monomial, the integration of $v$ with respect to a measure $\mu$ results in a moment of $\mu$. Therefore, \eqref{eq:LiouvillePolynomial} is a linear constraint on the moments of $\mu_0,$ $\mu$ and $\mu_T$. In this case, instead of finding a tuple of measures satisfying the constraints in~\eqref{eq:MainProgram}, we aim to find (finite) sequences of numbers that satisfy the constraint~\eqref{eq:LiouvillePolynomial}. Moreover, the sequences of numbers are moments of measures $\widetilde{\mu},\widehat{\mu},\mu,\mu_T.$ As required by~\eqref{eq:MainProgram}, these measures must be supported on certain specified sets. To formalize this idea, in order to obtain an approximated solution to~\eqref{eq:probsafe0}, we want to find sequences of numbers that are moments of the tuple of measures feasible in~\eqref{eq:MainProgram}. To better explain this approach, we first introduce necessary notions related to the multi-dimensional moment problem characterizing the relationship between sequences of numbers and moments of measures.

\subsection{Multi-dimensional $K$-moment problem}
Given an $\mathbb{R}^n$-valued random variable $\mathbf{x} \sim \nu$ and an integer vector $\boldsymbol{\alpha}\in \mathbb{N}^n,$ the $\boldsymbol{\alpha}$-moment of $\mathbf{x}$ is defined as $\mathbb{E}[\mathbf{x}^{\boldsymbol{\alpha}}]= \int_{\mathbb{R}^n}\prod_{i = 1}^n x_i^{\alpha_i}d\nu.$ Moreover, we define the \emph{order} of an $\boldsymbol{\alpha}$-moment to be $|\boldsymbol{\alpha}|$. Finally, a sequence $\mathbf{y}=\{y_{\boldsymbol{\alpha}}\}_{\boldsymbol{\alpha}\in \mathbb{N}^n}$ indexed by $\boldsymbol{\alpha}$ is called a \emph{multi-sequence}. Given a multi-sequence $\mathbf{y} = \{y_{\boldsymbol{\alpha}}\}_{\boldsymbol{\alpha}\in \mathbb{N}^n},$ we define the linear functional $L_\mathbf{y}: \mathbb{R}[\mathbf{x}] \rightarrow \mathbb{R}$ as
\begin{equation}\label{eq:Functional}
f(\mathbf{x}) = \sum_{\bdA \in \mathbb{N}^n} f_{\bdA}\mathbf{x}^{\bdA} \mapsto L_\mathbf{y}(f) = \sum_{\bdA \in \mathbb{N}^n}f_{\bdA}y_{\bdA}.
\end{equation}
The introduction of the above functional, often known as the \emph{Riesz functional}~\cite{lasserre2015introduction}, is convenient to express the moments of random variables. More specifically, let $\mathbf{x}$ be an $\mathbb{R}^n$-valued random variable with corresponding probability measure $\nu$ and let $f$ be a polynomial in $\mathbf{x}$. Then, the expectation of $f(\mathbf{x})$ is equal to
\begin{equation*}
\int f(\mathbf{x}) d\nu = \int \sum_{\bdA \in \mathbb{N}^n} f_{\bdA}\mathbf{x}^{\bdA}d\nu = \sum_{\bdA \in \mathbb{N}^n} f_{\bdA}y_{\bdA} = L_{\mathbf{y}}(f)
\end{equation*}
where $y_{\bdA}$ is the $\bdA$-moment of $\mathbf{x}.$

\begin{define}\label{def}
Let $K\subseteq \mathbb{R}^n$ be a closed set. Let $\mathbf{y} = \{y_{\boldsymbol{\alpha}}\}_{\boldsymbol{\alpha} \in \mathbb{N}^n}$ be an infinite real multi-sequence. A measure $\nu$ on $\mathbb{R}^n$ is said to be a $K$-representing measure for $\mathbf{y}$ if
\begin{equation}
y_{\boldsymbol{\alpha}} = \int_{\mathbb{R}^n} \mathbf{x}^{\boldsymbol{\alpha}}d\nu \text{ for all }\boldsymbol{\alpha}\in \mathbb{N}^n
\end{equation} 
and $\mathtt{supp}(\nu)\subseteq K.$ If $\mathbf{y}$ has a $K$-representing measure, we say that $\mathbf{y}$ is $K$-feasible. 
\end{define}
Note that not all multi-sequences are $K$-feasible, since there may not exist a measure supported on $K$ whose moments match the values in the multi-sequence. A necessary and sufficient condition for the feasibility of the $K$-moment problem, restricted to the case when $K$ is semi-algebraic and compact, can be stated in terms of linear matrix inequalities. These conditions involve \emph{moment matrices} and \emph{localizing matrices}, defined below. 
\begin{define}{\emph{\cite{bernard2009moments}}}\label{def:momentmatrix}
Let $\mathbf{y}_{n, 2r} = \{y_{\boldsymbol{\alpha}}\}_{\boldsymbol{\alpha} \in \mathbb{N}^n_{2r}}$ be a (finite) real multi-sequence. The \emph{moment matrix} of $\mathbf{y}_{n, 2r}$, denoted by $M_r(\mathbf{y}_{n, 2r})$, is defined as the real matrix indexed by $\mathbb{N}^n_r$ whose entries are
\begin{equation}\label{eq:MomentEntryWise}
[M_r(\mathbf{y}_{n, 2r})]_{\boldsymbol{\alpha}, \boldsymbol{\beta}} = y_{\boldsymbol{\alpha}+ \boldsymbol{\beta}}
\end{equation}
for all $\boldsymbol{\alpha}, \boldsymbol{\beta} \in \mathbb{N}^n_r.$
\end{define}
To better explain how the moment matrix is constructed, we consider $n = 2$, $r = 1$ and $\mathbf{y}_{2,2} = \{y_{00},y_{01},y_{10},y_{11},y_{02},y_{20} \}$ as an example. According to Definition~\ref{def:momentmatrix}, we have that
$$M_1(\mathbf{y}_{2,2}) = \begin{bmatrix}
y_{00} & y_{10} & y_{01} \\
y_{10} & y_{20} & y_{11}\\
y_{01} & y_{11} & y_{02}\\
\end{bmatrix}.$$ 



Similarly, we define the \emph{localizing matrices} as follows. 
\begin{define}\label{def:Localizing}
Consider a polynomial $g(\mathbf{x}) = \sum_{\boldsymbol{\gamma} \in \mathbb{N}^n} u_{\boldsymbol{\gamma}}\mathbf{x}^{\boldsymbol{\gamma}}.$ Given a finite multi-sequence $\mathbf{y}_{n, 2r} = \{y_{\boldsymbol{\alpha}}\}_{\boldsymbol{\alpha} \in \mathbb{N}^n_{2r}},$ the localizing matrix of $\mathbf{y}_{n, 2r}$ with respect to $g,$ denoted by $M_r(g,\mathbf{y}_{n, 2r}),$ is the real matrix indexed by $\mathbb{N}^n_r$ whose entries are
\begin{equation}
[M_r(g,\mathbf{y}_{n, 2r})]_{\boldsymbol{\alpha}, \boldsymbol{\beta}} = 
\sum_{\boldsymbol{\gamma} \in \mathbb{N}^n} u_{\boldsymbol{\gamma}} y_{\boldsymbol{\gamma} + \boldsymbol{\alpha} + \boldsymbol{\beta}}
\end{equation}
for all $\boldsymbol{\alpha},\boldsymbol{\beta} \in \mathbb{N}^n_r.$
\end{define}

Under specific assumptions on the set $K,$ it is possible to state necessary and sufficient conditions for $K$-feasibility of $\mathbf{y}$ using moment and localizing matrices. Such a method is built upon an algebraic characterization of the relationship between polynomials and sum-of-squares (SOS) polynomials.
\begin{define}{\emph{(Sum-of-squares polynomial)}}
A polynomial $p:\mathbb{R}^n \rightarrow \mathbb{R}$ is a sum-of-squares polynomial if $p$ can be written as 
\begin{equation}
p(\mathbf{x}) = \sum_{j \in J}p_j(\mathbf{x})^2, \ \mathbf{x} \in \mathbb{R}^n
\end{equation}
for some finite family of polynomials $\{p_j \mid j\in J\}.$
\end{define}
The following result utilizes the properties of sum-of-squares polynomials to characterize when a multi-sequence $\mathbf{y}$ is $K$-feasible in terms of moment and localizing matrices.
\begin{thm}\label{thm:PutinarCondition}{\emph{(Putinar's Positivstellensatz~\cite{putinar1993positive})}} Consider an infinite multi-sequence $\mathbf{y} = \{y_{\boldsymbol{\alpha}}\}_{\boldsymbol{\alpha} \in \mathbb{N}^n}$ and a collection of polynomials $g_i: \mathbb{R}^n \rightarrow \mathbb{R} \text{ for all } i \in [m].$ Define a compact semi-algebraic set $K = \left\{\mathbf{x}\in \mathbb{R}^n \mid g_i(\mathbf{x}) \ge 0, \ i\in[m]\right\}.$ Assume that there exists a polynomial $u = u_0 + \sum_{j=1}^m u_ig_i$ where $u_i$ are SOS polynomials for all $i \in \{ 0 \}  \cup [m]$ such that the set $\left\{\mathbf{x} \mid u(\mathbf{x}) \ge 0\right\}$ is compact. Then, $\mathbf{y}$ has a $K$-representing measure if and only if
\begin{equation}\label{eq:FeasibleConditionReduction}
\begin{aligned}
& M_r(\mathbf{y}) \succeq 0 \text{ and }\\
& M_r(g_{j},\mathbf{y}) \succeq 0, \text{ for all } j \in [m] \text{ and } r\in \mathbb{N}.
\end{aligned}
\end{equation}
\end{thm}
In the following subsection, we will leverage this theorem to construct approximate solutions of $\mathtt{P}$ and $\mathtt{D}$.


\subsection{Finite-dimensional approximations}\label{subsec:approximate}
\subsubsection{SDP relaxation of $\mathtt{P}$}
As mentioned above, in the relaxed version of $\mathtt{P}$, we aim to optimize over sequences of moments of a tuple of measures $(\widetilde{\mu}, \widehat{\mu},\mu,\mu_T)$. We use $(\widetilde{\mathbf{y}},\widehat{\mathbf{y}},\mathbf{y}, \mathbf{y}_T)$ to denote the moment sequences of the corresponding measures, respectively. On the one hand, since $\mu$ is supported on $[0, T]\times\mathcal{X},$ the elements in the moment sequence $\mathbf{y}$ are of the form $y_{\bdA}$ where $\bdA \in \mathbb{N} \times \mathbb{N}^n.$ On the other hand, since $\mu_T$ is supported on $\mathcal{X},$ the elements in $\mathbf{y}_T$ are of the form $y_{\bdA}$ where $\bdA \in \mathbb{N}^n.$ Using the Riesz functional~\eqref{eq:Functional} on~\eqref{eq:LiouvillePolynomial}, we obtain
\begin{equation}\label{eq:ConstraintLiouville}
\begin{aligned}
& L_{\mathbf{y}_T}(v(T, \cdot)) - L_{\mathbf{y}}(\mathcal{L}v) = L_{\mathbf{y}_0}(v(0, \cdot)) \\
& \text{for all } v(t, \mathbf{x}) = t^a \mathbf{x}^{\bdA} \text{ and } a+|\bdA| \le 2r.
\end{aligned}
\end{equation}
Applying the Riesz functional on the first linear constraint in $\mathtt{P},$ we have that
\begin{equation}\label{eq:ConstraintMeasureRelation}
\begin{aligned}
& L_{\widetilde{\mathbf{y}}}(w) + L_{\widehat{\mathbf{y}}}(w) = L_{\mathbf{y}}(w) \\
& \text{for all } w(t, \mathbf{x}) = t^a \mathbf{x}^{\bdA} \text{ and } a+|\bdA| \le 2r.
\end{aligned}
\end{equation}
Both equations in \eqref{eq:ConstraintMeasureRelation} are linear with respect to the elements in $\mathbf{y}, \widetilde{\mathbf{y}}, \widehat{\mathbf{y}}, \mathbf{y}_T$; hence, it is possible to write them compactly into a linear equation, as follows:
\begin{equation}
A_r(\widetilde{\mathbf{y}}, \widehat{\mathbf{y}}, \mathbf{y}, \mathbf{y}_T) = b_r.
\end{equation}

From Theorem~\ref{thm:PutinarCondition}, since $\mathtt{supp}(\mu) \subseteq [0, T]\times\mathcal{X},$ the moment and localizing matrices of $\mathbf{y}$ with respect to $g_i^{\mathcal{X}}$ are positive semidefinite for all positive integers $r \in \mathbb{N}.$ Let 
\begin{equation*}
	d_i^{\mathcal{X}_u} = \ceil{\frac{\text{deg} \ g_i^{\mathcal{X}_u}}{2}} \ \forall i \in [n_{\mathcal{X}_u}], \ d_j^{\mathcal{X}} = \ceil{\frac{\text{deg} \ g_j^{\mathcal{X}}}{2}} \ \forall j \in [n_{\mathcal{X}}]
\end{equation*}
where deg denotes the degree of a polynomial. Given a fixed positive integer $r \in\mathbb{N},$ we construct the $r$-th order relaxation of $\mathtt{P}$, as follows:
\begin{equation}\label{eq:PrimalSDP}
\begin{aligned}
\mathtt{P_{r}}: \indent  & \underset{( \widetilde{\mathbf{y}},\widehat{\mathbf{y}},\mathbf{y},\mathbf{y}_T)}{\text{\upshape{maximize}}} \ L_{\widetilde{\mathbf{y}}}(g)\\
& \text{\upshape{subject to}} \ A_r(\widetilde{\mathbf{y}}, \widehat{\mathbf{y}}, \mathbf{y}, \mathbf{y}_T) = b_r\\
& \mbox{\indent \indent \indent \indent} \ M_r(\widetilde{\mathbf{y}}) \succeq 0,\  M_{r-1}(t(T-t), \widetilde{\mathbf{y}}) \succeq 0 \\
& \mbox{\indent \indent \indent \indent} \ M_{r-d_i^{\mathcal{X}_u}}(g_i^{\mathcal{X}_u}, \widetilde{\mathbf{y}}) \succeq 0, \forall i \in [n_{\mathcal{X}_u}] \\
& \mbox{\indent \indent \indent \indent} \ M_r(\widehat{\mathbf{y}}) \succeq 0, M_{r-1}(t(T-t), \widehat{\mathbf{y}}) \succeq 0\\
& \mbox{\indent \indent \indent \indent} \ M_{r-d_i^{\mathcal{X}}}(g_i^{\mathcal{X}}, \widehat{\mathbf{y}}) \succeq 0, \forall i \in [n_{\mathcal{X}}] \\
& \mbox{\indent \indent \indent \indent} \ M_r(\mathbf{y}) \succeq 0, \ M_{r-1}(t(T-t), \mathbf{y}) \succeq 0\\
& \mbox{\indent \indent \indent \indent} \ M_{r-d_i^{\mathcal{X}}}(g_i^\mathcal{X}, \mathbf{y}) \succeq 0, \forall i \in [n_{\mathcal{X}}]\\
& \mbox{\indent \indent \indent \indent} \ M_r(\mathbf{y}_T) \succeq 0,\\
& \mbox{\indent \indent \indent \indent} \ M_{r-d_i^{\mathcal{X}}}(g_i^\mathcal{X}, \mathbf{y}_T) \succeq 0, \forall i \in [n_{\mathcal{X}}].\\
\end{aligned}
\end{equation}
In this program, the decision variable is the 4-tuple of finite multi-sequences $(\mathbf{\widetilde{y}}, \mathbf{\widehat{y}},\mathbf{y},\mathbf{y}_T).$ Furthermore, $\mathtt{P_r}$ is an SDP and, thus, can be solved using off-the-shelf software. In addition to relaxing the primal LP $\mathtt{P},$ it is also possible to relax the dual LP $\mathtt{D}$, as shown next.

\subsubsection{SOS relaxation of $\mathtt{D}$}
To formulate the relaxed program of $\mathtt{D},$ we begin by considering the dual of $\mathtt{P_r}.$ Furthermore, as shown in~$\mathtt{D},$ the decision variables are $v(t,\mathbf{x})\in \mathcal{C}^1([0,T]\times\mathcal{X})$ and $w(t,\mathbf{x}) \in \mathcal{C}([0,T]\times\mathcal{X}).$ The relaxed program is obtained by restricting the functions in \eqref{eq:DualProgram} to polynomials of degrees up to $2r$, and then replacing the non-negativity constraint with sum-of-squares constraints \cite{parrilo2000structured}. To formalize this argument, we first need to introduce some notations.

Given a semi-algebraic set $A = \lbrace \mathbf{x} \in \mathbb{R}^n \mid h_i(\mathbf{x}) \ge 0, h_i \in \mathbb{R}[\mathbf{x}], \forall i \in [m] \rbrace$, we define the $r$-th order quadratic module of $A$ as
\begin{equation}
\begin{aligned}
Q_r(A) = \big \lbrace q \in & \mathbb{R}[\mathbf{x}]_r \mid \exists \ \text{SOS} \ \lbrace s_k \rbrace_{k \in [m] \cup \lbrace 0 \rbrace} \subset \mathbb{R}[x]_{r} \\ & \text{s.t.} \ q = s_0 + \sum_{k \in [m]} h_k s_k \big \rbrace.
\end{aligned}
\end{equation}
Following a process similar to~\cite{mohan2016convex}, the relaxed dual program, denoted by $\mathtt{D}_r$, can be written as follows
\begin{equation} \label{eq:DualSOS}
\begin{aligned}
	\mathtt{D_r}: \indent \indent & 
     \text{\upshape{minimize}} \ \int v(0, \cdot) d\mu_0 \\
     & \text{\upshape{subject to}} \ w - g \in Q_{2r}([0,T] \times \mathcal{X}_u) \\
     & \mbox{\indent \indent \indent \indent} \ - \mathcal{L}v - w \in  Q_{2r}([0,T] \times \mathcal{X}) \\
     & \mbox{\indent \indent \indent \indent} \ v(T, \cdot) \in  Q_{2r}(\mathcal{X}) \\
     & \mbox{\indent \indent \indent \indent} \ w \in  Q_{2r}([0,T] \times \mathcal{X}).
\end{aligned}
\end{equation}
In this program, we optimize over the vector of polynomials $(w,v) \in \mathbb{R}[t,\mathbf{x}]_{2r} \times \mathbb{R}[t, \mathbf{x}]_{2r}$. 

Notice that $\mathtt{P}_r$ and $\mathtt{D}_r$ provide approximate solutions to $\mathtt{P}$ and $\mathtt{D},$ respectively. In the next theorem, we show that there is no duality gap between $\mathtt{P}_r$ and $\mathtt{D_r}$ and that the optimal values of $\mathtt{P}_r$ and $\mathtt{D}_r$ converge to the optimal values of $\mathtt{P}$ and $\mathtt{D}$, respectively, as $r$ increases.

\begin{thm}\label{thm:StrongDualityRelaxed}
Given a positive integer $r\in\mathbb{N},$ let $p_r^\star$ and $d_r^\star$ be the optimal values of $\mathtt{P}_r$ and $\mathtt{D}_r$, respectively. If $\mathcal{X}_u$ and $\mathcal{X}$ have nonempty interior, then $p_r^\star = d_r^\star.$ Furthermore, 
\begin{equation}
d_r^* = p_r^\star \downarrow p^\star = d^\star.
\end{equation}
\end{thm}
\begin{proof}
See appendix~\ref{appen:proof}.
\end{proof}
As a result of this theorem, $p^\star_r$ is a non-increasing function of $r$ and it converges asymptotically to $p^\star.$ From Theorem~\ref{thm:LPOptimal}, $p^\star$ is equal to the expected time the system spends in the unsafe region, as expressed in \eqref{eq:probsafe0}.

\section{Numerical Examples}
\begin{figure} 
\centering
	\includegraphics[width= 0.4\textwidth]{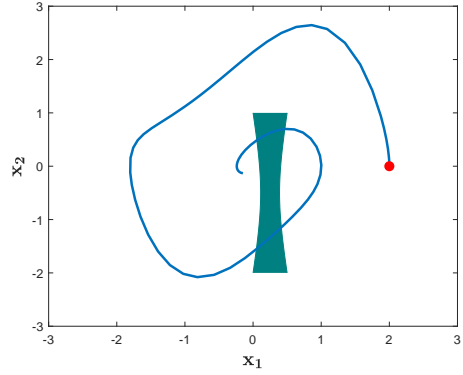}
	\caption{\label{fig:vanderpol} Trajectory $\mathbf{x}(t),$ where $t \in [0, 10],$ of the Van der Pol system (blue curve) with initial condition $\mathbf{x}_0 = [2, 0]^T$ (red circle). The unsafe region $\mathcal{X}_u$ is depicted by the nonconvex colored set.}
\end{figure}

\begin{figure} 
\centering
	\includegraphics[width= 0.4\textwidth]{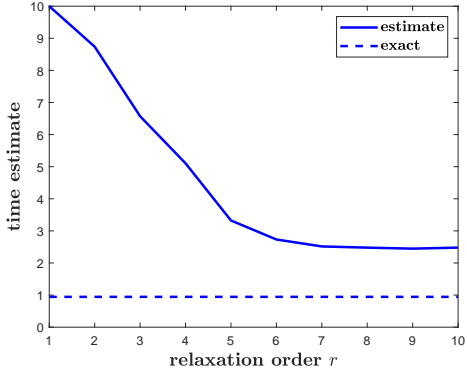}
	\caption{\label{fig:vanderpol_sol} This figure shows the exact value (dashed line) and the approximation (solid line) to~\eqref{eq:probsafe0} using $\mathtt{D}_r$ with different values of $r$. The system dynamics under consideration is the Van der Pol system~\eqref{eq:vanderpolsys}, whereas the initial distribution is $\mu_0 = \delta_{[2, 0]^\top}.$}
\end{figure}

In this section, we provide a numerical example to illustrate our framework. We complete all numerical simulations using YALMIP \cite{lofberg2004yalmip} (for sum-of-squares programs) and MOSEK \cite{mosek2015mosek} (for semidefinite programs). 
In particular, we evaluate our framework on the Van der Pol oscillator -- a second order nonlinear dynamical system whose dynamics is given by
\begin{equation} \label{eq:vanderpolsys}
	\begin{aligned}
	&\dot{x}_1 = - x_2 \\
	&\dot{x}_2 = x_1 + (x_1^2 - 1)x_2.
	\end{aligned}
\end{equation}
Moreover, we consider the following parameter settings (see Figure~\ref{fig:vanderpol}): (i) the final time is set to be $T = 10$, (ii) the initial condition is set to be $\mathbf{x}(0) = \mathbf{x}_0 = [2, 0]^\top$, and (iii) the unsafe region is specified by a nonconvex two-dimensional semi-algebraic set $\mathcal{X}_u = \{(x_1, x_2) \in \mathbb{R}^2 \mid 52(x_1 - 0.25)^2 - (x_2 + 0.5)^2 \leq 1, 0 \leq x_1 \leq 0.5, -2 \leq x_2 \leq 1 \}$. To ease the numerical computations, we adopt proper scaling of the system's coordinates such that $T$ and $\mathcal{X}$ are normalized to be $T = 1$ and $\mathcal{X} = [-1, 1] \times [-1, 1]$, respectively. In this case,~\eqref{eq:probsafe0} cannot be computed analytically. However, through numerical simulation, we obtain that the Van der Pol oscillator spends (approximately) $0.9446$ seconds in the unsafe region $\mathcal{X}_u$. We demonstrate our upper bounds on this time using $\mathtt{D}_r$ with varying values of $r$ in Figure~\ref{fig:vanderpol_sol}.


\section{Conclusion}
In this paper, we have proposed a flexible safety verification notion for nonlinear autonomous systems described via polynomial dynamics and unsafe regions described via polynomial inequalities. Instead of verifying safety by checking whether the dynamics completely avoids the unsafe regions, we consider the system to be safe if it spends less than a certain amount of time in these regions. This more flexible notion can be of relevance in, for example, solar-powered vehicles where the vehicle should avoid spending too much time is dark areas. More generally, this framework can be useful in those situations where the system is able to tolerate the exposure to a deteriorating agent, such as excessive heat or radiation, for a limited amount of time. In this paper, we first propose an infinite-dimensional LP over the space of measures whose solution is equal to the (expected) time our (nonlinear) system spends in the (possibly nonconvex) unsafe regions. We then approximate the solution of the LP through a monotonically converging sequence of upper bounds by solving a hierarchy of SDPs. We have validated our approach via a simple example involving a nonlinear Van der Pol oscillator. As future work, we are working on the problem of path planning using the flexible safety notion herein proposed.

\begin{appendices}
\section{}\label{appen:proof}
\begin{proof}[Proof of Theorem~\ref{thm:LPOptimal}]
	First, we show that when the initial distribution $\mu_0$ and the system dynamics \eqref{eq:dynamics} are given, the Liouville equation \eqref{eq:AverageLiouvilleLinear} has a unique solution $(\mu, \mu_T)$ up to a subset of $[0,T] \times \mathcal{X}$ of Lebesgue measure zero and $(\mu, \mu_T)$ coincide with the average occupation measure defined by \eqref{eq:AverageOccupationMeasure} and the average final measure defined by \eqref{eq:AverageFinalMeasure}. Let $(\mu, \mu_T)$ be a pair of measures satisfying \eqref{eq:AverageLiouvilleLinear}. From \cite[Lemma~3]{henrion2014convex}, $\mu$ can be disintegrated as $d\mu(t,\mathbf{x}) = d\mu_t(\mathbf{x})dt$ where $dt$ is the Lebesgue measure on $[0,T]$. $\mu_t(\mathbf{x})$ is a stochastic kernel on $\mathcal{X}$ given $t$ and can be interpreted as the distribution of the states at time $t$ following the evolution of \eqref{eq:dynamics} with $\mathbf{x}_0 \sim \mu_0$. $\mu_t(\mathbf{x})$ is uniquely defined $dt$-almost everywhere. As proved in \cite[Lemma~3]{henrion2014convex}, $\mu_t$ satisfies a continuity equation which implies $\mu$ and $\mu_T$ coincide with the average occupation measure and the average final measure generated by the family of absolutely continuous admissible trajectories of \eqref{eq:dynamics} starting from $\mu_0$.
	
	Then solving $\mathtt{P}$ can be decomposed into two steps: first find a feasible $(\mu, \mu_T) \in \mathcal{M}_{+}([0,T] \times \mathcal{X}) \times \mathcal{M}_{+}(\mathcal{X})$ to the Liouville equation $\delta_T \otimes \mu_T = \delta_0 \otimes \mu_0 + \mathcal{L}^*\mu$ and then solve the following optimization problem:
	\begin{equation}\label{eq:VolumeProblem}
\mathtt{Q}: \indent	\underset{\widetilde{\mu}}{\sup} \lbrace \int g d\widetilde{\mu} : \widetilde{\mu} \leq \mu; \widetilde{\mu} \in \mathcal{M}_{+}([0,T] \times \mathcal{X}_u) \rbrace.
	\end{equation}
	Since $\mathcal{X}$ and $\mathcal{X}_u$ are compact with $\mathcal{X}_u \subseteq \mathcal{X}$, by \cite[Theorem~3.1]{henrion2009approximate} the restriction $\widetilde{\mu}^*$ of $\mu$ to $\mathcal{X}_u$ defined by \eqref{eq:VolumeDefinition} is the unique optimal solution to $\mathtt{Q}$ and $\sup \mathtt{Q} = \max \mathtt{Q} = \int g d\widetilde{\mu}^* = \int_{\mathcal{X}_u} g d\mu$. 
	
	As the feasible $\mu$ in $\mathtt{P}$ coincides with the average occupation measure in \eqref{eq:AverageOccupationMeasure}, $\widetilde{\mu}^*$ is also the $\widetilde{\mu}$-component of an optimal solution to $\mathtt{P}$ and $\sup \mathtt{P} = \max \mathtt{P} = \int g d\widetilde{\mu}^*$. When $g \equiv 1$, we have $\max \mathtt{P} = \mu([0,T] \times \mathcal{X}_u)$ with $\mu$ being the average occupation measure defined in \eqref{eq:AverageOccupationMeasure}. 
\end{proof}

\begin{proof}[Proof of Theorem~\ref{thm:StrongDuality}]
	The proof follows the same lines as that of \cite[Theorem~2]{henrion2014convex}. Define 
	\begin{equation*}
		\mathbf{C} = \mathcal{C}([0,T] \times \mathcal{X}_u) \times \mathcal{C}([0,T] \times \mathcal{X}) \times \mathcal{C}([0,T] \times \mathcal{X}) \times \mathcal{C}(\mathcal{X})
	\end{equation*}
	\begin{equation*}
		\mathbf{M} \!=\! \mathcal{M}([0,T] \times \!\mathcal{X}_u) \!\times \!\mathcal{M}([0,T] \times \mathcal{X}) \times \mathcal{M}([0,T] \times \mathcal{X}) \times \mathcal{M}(\mathcal{X})
	\end{equation*}
	and let $\mathcal{K}$ and $\mathcal{K}'$ denote the positive cones of $\mathbf{C}$ and $\mathbf{M}$, respectively. By Riesz-Markov-Kakutani representation theorem~\cite{kakutani1941concrete}, $\mathcal{K}'$ is the topological dual of the cone $\mathcal{K}$. The infinite dimensional linear program $\mathtt{P}$ can be written as:
	\begin{equation} \label{eq:PrimalProblem}
	\begin{aligned}
	\sup \indent & \langle \gamma, c \rangle \\
	\text{s.t.} \indent & \mathcal{A}' \gamma = \beta, \quad \gamma \in \mathcal{K}'
	\end{aligned}	
	\end{equation}
	where the supremum is taken over the vector $\gamma = (\widetilde{\mu}, \widehat{\mu}, \mu, \mu_T)$, the linear operator $\mathcal{A}': \mathcal{K}' \rightarrow \mathcal{C}^{1}([0,T] \times \mathcal{X})^* \times \mathcal{M}([0,T] \times \mathcal{X})$ is defined by $\mathcal{A}' \gamma = (\delta_T \otimes \mu_T - \mathcal{L}^*\mu, \mu - \widetilde{\mu} - \widehat{\mu})$ and $\beta = (\delta_0 \otimes \mu_0, 0) \in \mathcal{C}^{1}([0,T] \times \mathcal{X})^* \times \mathcal{M}([0,T] \times \mathcal{X}) $. The vector of functions in the objective is $c = (g, 0, 0, 0)$. Define the duality bracket between a vector of measures $\nu \in (\mathcal{M}(\mathcal{S}))^p$ and a vector of functions $h \in (\mathcal{C}(\mathcal{S}))^p$ over a topological space $\mathcal{S}$ by $\langle h, \nu \rangle = \sum_{i = 1}^p \int_{\mathcal{S}} {[h]}_i d{[\nu]}_i$. Then $\langle \gamma, c \rangle = \int g d\widetilde{\mu}$.
	
	The dual to \eqref{eq:PrimalProblem} can be interpreted as:
	\begin{equation} \label{eq:DualProblem}
	\begin{aligned}
	\inf \indent & \langle \beta, z \rangle \\
	\text{s.t.} \indent & \mathcal{A}z - c \in \mathcal{K}
	\end{aligned}
	\end{equation} 
	where the infimum is over $z = (v, w) \in \mathcal{C}^{1}([0,T] \times \mathcal{X}) \times \mathcal{C}([0,T] \times \mathcal{X})$, the linear operator $\mathcal{A}:  \mathcal{C}^{1}([0,T] \times \mathcal{X}) \times \mathcal{C}([0,T] \times \mathcal{X}) \rightarrow \mathbf{C}$ is given by $\mathcal{A}z = (w, w, - \mathcal{L}v - w, v(T, \cdot))$ and satisfies the adjoint property $\langle \mathcal{A}' \gamma, z \rangle = \langle \gamma, \mathcal{A}z \rangle$. The linear program \eqref{eq:DualProblem} is exactly \eqref{eq:DualProgram}.
	
	From \cite[Theorem~3.10]{anderson1987linear}, there is no duality gap between \eqref{eq:PrimalProblem} and \eqref{eq:DualProblem} if the supremum of \eqref{eq:PrimalProblem} is finite and the set $P = \lbrace (\mathcal{A}' \gamma, \langle \gamma, c \rangle) \mid \gamma \in \mathcal{K}' \rbrace$ is closed in the weak* topology of $\mathcal{K}'$. Since $\widetilde{\mu}$ is dominated by the average occupation measure $\mu$ and its underlying support is compact, the supremum of \eqref{eq:PrimalProblem} is finite. To prove that $P$ is closed, consider a sequence $\gamma_k = (\widetilde{\mu}^k, \widehat{\mu}^k, \mu^k, \mu_T^k) \in \mathcal{K}'$ such that $\mathcal{A}' \gamma_k \rightarrow a$ and $\langle \gamma_k, c \rangle \rightarrow b$ as $k \rightarrow \infty$ for some $(a,b) \in \mathcal{C}^{1}([0,T] \times \mathcal{X})^* \times \mathcal{M}([0,T] \times \mathcal{X}) \times \mathbb{R}$. Consider the test function $z_1 = (T - t, 0)$ which gives $\langle \mathcal{A}' \gamma_k, z_1 \rangle = \mu^k([0,T] \times \mathcal{X}) \rightarrow \langle a,z_1 \rangle < \infty$; since the measures $\mu^k$ are non-negative, this implies $\{ \mu^k \}$ is bounded. By taking $z_2 = (1,-1)$, we have $\langle \mathcal{A}' \gamma_k, z_2 \rangle = \mu_T^k(\mathcal{X}) + \widetilde{\mu}^k([0,T] \times \mathcal{X}_u) + \widehat{\mu}^k([0,T] \times \mathcal{X}) - \mu^k([0,T] \times \mathcal{X}) \rightarrow \langle a,z_2 \rangle < \infty$; since $\{ \mu^k \}$ is bounded, by similar arguments the sequences $\{ \widetilde{\mu}^k \}$, $\{ \widehat{\mu}^k  \}$ and $\{ \mu_T^k \}$ are bounded as well.
	
	As a result, $\lbrace \gamma_k \rbrace$ is bounded and we can find a ball $B$ in $\mathbf{M}$ with $\lbrace \gamma_k \rbrace \subset B$. From the weak* compactness of the unit ball (Alaoglu's theorem \cite[Section~5.10, Theorem~1]{luenberger1997optimization}) there is a subsequence $\lbrace \gamma_{k_i} \rbrace$ that weak*-converges to some $\gamma \in \mathcal{K}'$. Notice that $\mathcal{A}'$ is weak*-continuous because $\mathcal{A}z \in \mathbf{C}$ for all $z \in \mathcal{C}^{1}([0,T] \times \mathcal{X}) \times \mathcal{C}([0,T] \times \mathcal{X})$. So $(a, b) = \lim_{i \rightarrow \infty} (\mathcal{A}' \gamma_{k_i}, \langle \gamma_{k_i}, c \rangle) = (\mathcal{A}' \gamma, \langle \gamma, c \rangle) \in P$ by the continuity of $\mathcal{A}'$ and $P$ is closed.
\end{proof}

\begin{proof}[Proof of Theorem~\ref{thm:StrongDualityRelaxed}]
	The proof of strong duality follows from standard SDP duality theory. Let $\Delta_{\mu} = (\widetilde{\mu}, \widehat{\mu}, \mu, \mu_T)$ be the optimal solution to $\mathtt{P}$ and $\Delta_y = (\mathbf{\widetilde{y}}, \mathbf{\widehat{y}}, \mathbf{y}, \mathbf{y_T})$ be their corresponding moment sequences. Any finite truncation of $\Delta_y$ gives a feasible solution to $\mathtt{P_r}$. As $\mathcal{X}$ and $\mathcal{X}_u$ have non-empty interior, we have the truncation of $\Delta_y$ is strictly feasible for $\mathtt{P_r}$. By Slater's condition \cite{boyd2004convex}, there is no duality gap between $\mathtt{P_r}$ and $\mathtt{D_r}$, i.e., $p_r^* = d_r^*$. 
	
	The proof of convergence follows from \cite[Theorem~3.6]{lasserre2008nonlinear}. Since $[0,T]$, $\mathcal{X}$ and $\mathcal{X}_u$ are compact sets, we can assume after appropriate scaling $T = 1$ and $\mathcal{X} \times \mathcal{X}_u \subseteq [-1,1]^{n_{\mathcal{X}}} \times [-1,1]^{n_{\mathcal{X}_u}}$, which implies that the feasible set of the semidefinite program $\mathtt{P_r}$ is compact. Let $\Delta_r^* = (\mathbf{\widetilde{y}}_r^*, \mathbf{\widehat{y}}_r^*, \mathbf{y}_r^*, \mathbf{y_T}_r^*)$ be the optimal solution of $\mathtt{P_r}$ and complete the finite vectors $(\mathbf{\widetilde{y}}_r^*, \mathbf{\widehat{y}}_r^*, \mathbf{y}_r^*, \mathbf{y_T}_r^*)$ with zeros to make them infinite sequences. By a standard diagonal argument, there is a subsequence $\lbrace r_k \rbrace$ and a tuple of infinite vectors $\Delta^* = (\mathbf{\widetilde{y}}^*, \mathbf{\widehat{y}}^*, \mathbf{y}^*, \mathbf{y_T}^*)$ such that $\Delta_{r_k}^* \rightarrow \Delta^*$ as $k \rightarrow \infty$, where the convergence is interpreted as elementary-wise. Since the infinite vector $\mathbf{\widetilde{y}}^*$ in $\Delta^*$ is the limit point of a subsequence of the optimal solutions $\mathbf{\widetilde{y}}_r^*$ of $\mathtt{P_r}$, $\mathbf{\widetilde{y}}^*$ satisfies all the constraints in $\mathtt{P_r}$ as $r \rightarrow \infty$. Then by \emph{Putinar's Positivstellensatz}, $\mathbf{\widetilde{y}}^*$ has a representing measure $\widetilde{\mu}^*$ supported on $[0,T] \times \mathcal{X}_u$. Similarly, $\mathbf{\widehat{y}}^*, \mathbf{y}^*$ and $\mathbf{y_T}^*$ have their representing measures $\widehat{\mu}^*, \mu^*$ and $\mu_T^*$ with corresponding supports, respectively.
	
    As problem $\mathtt{P_r}$ is a relaxation of $\mathtt{P}$, $p_r^* \geq p^*$ for each $r$. Thus we have $\lim_{k \rightarrow \infty} \sup \mathtt{P_{r_k}} = \lim_{k \rightarrow \infty} L_{\mathbf{\widetilde{y}}_{r_k}^*}(g) = L_{\mathbf{\widetilde{y}}^*}(g) = \int g d\widetilde{\mu}^* \geq p^*$. On the other hand, $\mathcal{A}_r(\Delta^*) = \lim_{k \rightarrow \infty} \mathcal{A}_r(\Delta_{r_k}^*) = b_r$ for each $r \in \mathbb{N}$. Let $(\widetilde{\mu}^*, \widehat{\mu}^*, \mu^*, \mu_T^*)$ be the tuple of representing measures of $\Delta^*$. As measures on compact sets are determined by moments, $(\widetilde{\mu}^*, \widehat{\mu}^*, \mu^*, \mu_T^*)$ is a feasible solution to $\mathtt{P}$ which implies $\int g d \widetilde{\mu}^* \leq p^*$. Hence $\int g d \widetilde{\mu}^* = p^*$
	and $(\widetilde{\mu}^*, \widehat{\mu}^*, \mu^*, \mu_T^*)$ is an optimal solution of $\mathtt{P}$. For any $r$ we have $p_r^* \geq p_{r+1}^*$ because as $r$ increases, the constraints in $\mathtt{P_r}$ become more restrict. As a result, $p_{r_k}^* \downarrow p^*$ and furthermore $p_r^* \downarrow p^*$. By strong duality, $d_r^* = p_r^* \downarrow p^* = d^*$.
\end{proof}
\end{appendices}

\small
\bibliographystyle{IEEEtran}
{ \small \bibliography{reference}}

\begin{thebibliography}{10}
\providecommand{\url}[1]{#1}
\csname url@samestyle\endcsname
\providecommand{\newblock}{\relax}
\providecommand{\bibinfo}[2]{#2}
\providecommand{\BIBentrySTDinterwordspacing}{\spaceskip=0pt\relax}
\providecommand{\BIBentryALTinterwordstretchfactor}{4}
\providecommand{\BIBentryALTinterwordspacing}{\spaceskip=\fontdimen2\font plus
\BIBentryALTinterwordstretchfactor\fontdimen3\font minus
  \fontdimen4\font\relax}
\providecommand{\BIBforeignlanguage}[2]{{%
\expandafter\ifx\csname l@#1\endcsname\relax
\typeout{** WARNING: IEEEtran.bst: No hyphenation pattern has been}%
\typeout{** loaded for the language `#1'. Using the pattern for}%
\typeout{** the default language instead.}%
\else
\language=\csname l@#1\endcsname
\fi
#2}}
\providecommand{\BIBdecl}{\relax}
\BIBdecl

\bibitem{hu2003probabilistic}
J.~Hu, M.~Prandini, and S.~Sastry, ``Probabilistic safety analysis in three
  dimensional aircraft flight,'' in \emph{Proceedings of IEEE Conference on
  Decision and Control}, vol.~5.\hskip 1em plus 0.5em minus 0.4em\relax IEEE,
  2003, pp. 5335--5340.

\bibitem{glavaski2005safety}
S.~Glavaski, A.~Papachristodoulou, and K.~Ariyur, ``Safety verification of
  controlled advanced life support system using barrier certificates,'' in
  \emph{International Workshop on Hybrid Systems: Computation and
  Control}.\hskip 1em plus 0.5em minus 0.4em\relax Springer, 2005, pp.
  306--321.

\bibitem{ziegler2010fast}
J.~Ziegler and C.~Stiller, ``Fast collision checking for intelligent vehicle
  motion planning,'' in \emph{Intelligent Vehicles Symposium (IV), 2010
  IEEE}.\hskip 1em plus 0.5em minus 0.4em\relax IEEE, 2010, pp. 518--522.

\bibitem{althoff2007safety}
M.~Althoff, O.~Stursberg, and M.~Buss, ``Safety assessment of autonomous cars
  using verification techniques,'' in \emph{Proceedings of American Control
  Conference}.\hskip 1em plus 0.5em minus 0.4em\relax IEEE, 2007, pp.
  4154--4159.

\bibitem{althoff2009model}
------, ``Model-based probabilistic collision detection in autonomous
  driving,'' \emph{IEEE Transactions on Intelligent Transportation Systems},
  vol.~10, no.~2, pp. 299--310, 2009.

\bibitem{bemporad2000optimization}
A.~Bemporad, F.~D. Torrisi, and M.~Morari, ``Optimization-based verification
  and stability characterization of piecewise affine and hybrid systems,'' in
  \emph{International Workshop on Hybrid Systems: Computation and
  Control}.\hskip 1em plus 0.5em minus 0.4em\relax Springer, 2000, pp. 45--58.

\bibitem{chutinan2003computational}
A.~Chutinan and B.~H. Krogh, ``Computational techniques for hybrid system
  verification,'' \emph{IEEE Transactions on Automatic Control}, vol.~48,
  no.~1, pp. 64--75, 2003.

\bibitem{bollobas1997volume}
B.Bollob{\'a}s, ``Volume estimates and rapid mixing,'' \emph{Flavors of
  geometry}, vol.~31, pp. 151--182, 1997.

\bibitem{dyer1988complexity}
M.~E. Dyer and A.~M. Frieze, ``On the complexity of computing the volume of a
  polyhedron,'' \emph{SIAM Journal on Computing}, vol.~17, no.~5, pp. 967--974,
  1988.

\bibitem{anai2001reach}
H.~Anai and V.~Weispfenning, ``Reach set computations using real quantifier
  elimination,'' in \emph{International Workshop on Hybrid Systems: Computation
  and Control}.\hskip 1em plus 0.5em minus 0.4em\relax Springer, 2001, pp.
  63--76.

\bibitem{tomlin2003computational}
C.~J. Tomlin, I.~Mitchell, A.~M. Bayen, and M.~Oishi, ``Computational
  techniques for the verification of hybrid systems,'' \emph{Proceedings of the
  IEEE}, vol.~91, no.~7, pp. 986--1001, 2003.

\bibitem{asarin2003reachability}
E.~Asarin, T.~Dang, and A.~Girard, ``Reachability analysis of nonlinear systems
  using conservative approximation,'' in \emph{International Workshop on Hybrid
  Systems: Computation and Control}.\hskip 1em plus 0.5em minus 0.4em\relax
  Springer, 2003, pp. 20--35.

\bibitem{prajna2004safety}
S.~Prajna and A.~Jadbabaie, ``Safety verification of hybrid systems using
  barrier certificates,'' in \emph{International Workshop on Hybrid Systems:
  Computation and Control}.\hskip 1em plus 0.5em minus 0.4em\relax Springer,
  2004, pp. 477--492.

\bibitem{prajna2006barrier}
S.~Prajna, ``Barrier certificates for nonlinear model validation,''
  \emph{Automatica}, vol.~42, no.~1, pp. 117--126, 2006.

\bibitem{prajna2007framework}
S.~Prajna, A.~Jadbabaie, and G.~J. Pappas, ``A framework for worst-case and
  stochastic safety verification using barrier certificates,'' \emph{IEEE
  Transactions on Automatic Control}, vol.~52, no.~8, pp. 1415--1428, 2007.

\bibitem{sloth2012compositional}
C.~Sloth, G.~J. Pappas, and R.~Wisniewski, ``Compositional safety analysis
  using barrier certificates,'' in \emph{Proceedings of the 15th ACM
  international conference on Hybrid Systems: Computation and Control}.\hskip
  1em plus 0.5em minus 0.4em\relax Citeseer, 2012, pp. 15--24.

\bibitem{kousik2017safe}
S.~Kousik, S.~Vaskov, M.~Johnson-Roberson, and R.~Vasudevan, ``Safe trajectory
  synthesis for autonomous driving in unforeseen environments,'' in \emph{ASME
  2017 Dynamic Systems and Control Conference}.\hskip 1em plus 0.5em minus
  0.4em\relax American Society of Mechanical Engineers, 2017.

\bibitem{vinter1993convex}
R.~Vinter, ``Convex duality and nonlinear optimal control,'' \emph{SIAM Journal
  on Control and Optimization}, vol.~31, no.~2, pp. 518--538, 1993.

\bibitem{henrion2009approximate}
D.~Henrion, J.~B. Lasserre, and C.~Savorgnan, ``Approximate volume and
  integration for basic semialgebraic sets,'' \emph{SIAM Review}, vol.~51,
  no.~4, pp. 722--743, 2009.

\bibitem{lasserre2008nonlinear}
J.~B. Lasserre, D.~Henrion, C.~Prieur, and E.~Tr{\'e}lat, ``Nonlinear optimal
  control via occupation measures and lmi-relaxations,'' \emph{SIAM Journal on
  Control and Optimization}, vol.~47, no.~4, pp. 1643--1666, 2008.

\bibitem{henrion2014convex}
D.~Henrion and M.~Korda, ``Convex computation of the region of attraction of
  polynomial control systems,'' \emph{IEEE Transactions on Automatic Control},
  vol.~59, no.~2, pp. 297--312, 2014.

\bibitem{majumdar2014convex}
A.~Majumdar, R.~Vasudevan, M.~M. Tobenkin, and R.~Tedrake, ``Convex
  optimization of nonlinear feedback controllers via occupation measures,''
  \emph{The International Journal of Robotics Research}, vol.~33, no.~9, pp.
  1209--1230, 2014.

\bibitem{kakutani1941concrete}
S.~Kakutani, ``Concrete representation of abstract (m)-spaces (a
  characterization of the space of continuous functions),'' \emph{Annals of
  Mathematics}, pp. 994--1024, 1941.

\bibitem{arnol2013mathematical}
V.~I. Arnol'd, \emph{Mathematical methods of classical mechanics}.\hskip 1em
  plus 0.5em minus 0.4em\relax Springer Science \& Business Media, 2013,
  vol.~60.

\bibitem{zhao2017control}
P.~Zhao, S.~Mohan, and R.~Vasudevan, ``Control synthesis for nonlinear optimal
  control via convex relaxations,'' in \emph{Proceedings of American Control
  Conference}.\hskip 1em plus 0.5em minus 0.4em\relax IEEE, 2017, pp.
  2654--2661.

\bibitem{bernard2009moments}
J.~B. Lasserre, \emph{Moments, positive polynomials and their
  applications}.\hskip 1em plus 0.5em minus 0.4em\relax World Scientific, 2009,
  vol.~1.

\bibitem{lasserre2015introduction}
------, \emph{An introduction to polynomial and semi-algebraic
  optimization}.\hskip 1em plus 0.5em minus 0.4em\relax Cambridge University
  Press, 2015, vol.~52.

\bibitem{putinar1993positive}
M.~Putinar, ``Positive polynomials on compact semi-algebraic sets,''
  \emph{Indiana University Mathematics Journal}, vol.~42, no.~3, pp. 969--984,
  1993.

\bibitem{parrilo2000structured}
P.~A. Parrilo, ``Structured semidefinite programs and semialgebraic geometry
  methods in robustness and optimization,'' Ph.D. dissertation, California
  Institute of Technology, 2000.

\bibitem{mohan2016convex}
S.~Mohan and R.~Vasudevan, ``Convex computation of the reachable set for hybrid
  systems with parametric uncertainty,'' in \emph{Proceedings of American
  Control Conference}.\hskip 1em plus 0.5em minus 0.4em\relax IEEE, 2016, pp.
  5141--5147.

\bibitem{lofberg2004yalmip}
J.~L{\"o}fberg, ``Yalmip: A toolbox for modeling and optimization in matlab,''
  in \emph{Proceedings of the CACSD Conference}, vol.~3.\hskip 1em plus 0.5em
  minus 0.4em\relax Taipei, Taiwan, 2004.

\bibitem{mosek2015mosek}
A.~Mosek, ``The mosek optimization toolbox for matlab manual,'' 2015.

\bibitem{anderson1987linear}
E.~J. Anderson and P.~Nash, \emph{Linear programming in infinite-dimensional
  spaces: theory and applications}.\hskip 1em plus 0.5em minus 0.4em\relax John
  Wiley \& Sons, 1987.

\bibitem{luenberger1997optimization}
D.~G. Luenberger, \emph{Optimization by vector space methods}.\hskip 1em plus
  0.5em minus 0.4em\relax John Wiley \& Sons, 1997.

\bibitem{boyd2004convex}
S.~Boyd and L.~Vandenberghe, \emph{Convex optimization}.\hskip 1em plus 0.5em
  minus 0.4em\relax Cambridge university press, 2004.

\end{thebibliography}

\end{document}